\documentclass[12pt]{article}
\begin{document}
%\renewcommand{\baselinestretch}{1.6}
%\preprint{}
%\begin{frontmatter}
\title{Twin Primes and the Zeros of the Riemann Zeta Function}
\author{H. J. Weber\\Department of Physics\\
University of Virginia\\Charlottesville, VA 22904\\USA}
%\date{\today}
\maketitle
\begin{abstract}
The Legendre-type relation for the counting function of ordinary 
twin primes is reworked in terms of the inverse of the Riemann 
zeta function. Its analysis sheds light on the distribution of 
the zeros of the Riemann zeta function in the critical strip 
and their link to the twin prime problem. 
\end{abstract}
\vspace{3ex}
%\pacs{}
\leftline{MSC: 11A41, 11N05, 11M06}
\leftline{Keywords: zeta zeros, twin primes, non-ranks, Perron integrals} 

%\vspace{3ex}

\section{Introduction}

The pair sieve for ordinary twin primes~\cite{adhjw} 
leads to a formula for the twin prime counting function 
$\pi_2(x)$ that is analogous to Legendre's formula~\cite{hri} 
for the prime number counting function $\pi(x).$ Before and 
after separating it into main and error terms~\cite{adhjw}, 
it is rewritten here using relevant Dirichlet series. Since 
the Riemann zeta function ends up in the denominator of the 
contour integrals, this feature links the zeta zeros to twin 
primes, much like $\pi(x)$ or related counting functions 
are expressed as Perron integrals over $\zeta'/\zeta$ in 
analytic number theory~\cite{rm},\cite{hw}. Our analysis 
sheds light on the role of twin primes in the distribution 
of the nontrivial zeros of the Riemann zeta function, which 
are those in the critical strip, as usual. 
 
In Sect.~2 the main concepts, such as twin ranks, non-ranks and 
remnants of the twin-prime pair sieve are recalled along with its 
main result, the Legendre type formula for $\pi_2.$ In Sect.~3 
it is rewritten as a Perron integral and analyzed. In Sect.~4 
the findings are summarized and discussed.   

\section{Review of the Pair Sieve and Notations}    

The prime numbers $2, 3$ do not play an active role here because 
they are not of the standard form $6m\pm 1.$ This also applies to  
the first twin prime pair $3, 5.$ From now on $p$ denotes a prime 
number or variable and $p_j$ the $j$th prime with $p_1=2, p_2=3, 
p_3=5,\ldots .$ In our twin prime sieve $p_j$ plays the role of 
the variable $\sqrt{x}$ in Eratosthenes' sieve. 

{\bf Definition~2.1.} If $6m\pm 1$ is an ordinary twin prime pair for 
some positive integer $m$, then $m$ is its {\it twin rank}. A positive 
integer $n$ is a {\it non-rank} if $6n\pm 1$ are {\bf not both} prime. 

The arithmetical function $N(x),~x\not =n+\frac{1}{2}$ is needed for 
non-ranks. 

{\bf Definition~2.2.} If $x$ is real then $N(x)$ is the integer 
nearest to $x.$ The ambiguity for $x=n+\frac{1}{2}$ with integral 
$n$ will not arise.    

In Ref.~\cite{adhjw} we then prove 

{\bf Lemma~2.3.} {\it If $p\geq 5$ is prime then the positive 
integers   
\begin{eqnarray}
k(n,p)^{\pm}=n p\pm N(\frac{p}{6})>0,~n=0, 1, 2,\ldots 
\label{nr}
\end{eqnarray}
are non-ranks. If an integer $k>0$ is a non-rank, there is a 
prime $p\geq 5$ so that Eq.~(\ref{nr}) holds with either $+$ 
or $-$ sign.}  

This means that the pairs $6k^+\pm 1$ and $6k^-\pm 1$ each contain 
at least one composite number. Therefore, the primes $p\geq 5$ 
organize all non-rank numbers in pairs of arithmetic progressions.  
These pairs are twin prime analogs of multiples $n p,~n>1,$ of 
primes struck from the integers in Eratosthenes' sieve.     

Given a prime $p\geq 5,$ when all non-ranks to primes $5\leq p'<p$ 
are subtracted from the non-ranks to $p,$ then the non-ranks to 
{\bf parent} prime $p$ are left forming the set ${\cal A}_p.$ 
This process~\cite{adhjw} naturally introduces the primorial 
$L(p)=\prod_{5\leq p'<p}p'$ as the period (of its arithmetic 
progressions). $L(p_j)\to\infty$ is the twin prime sieve's 
analog of the variable $x\to\infty$ in Eratosthenes' sieve.  

{\bf Definition~2.4.} Let $p\geq p'\geq 5$ be prime. The supergroup 
${\cal S}_p=\bigcup_{p'\leq p}{\cal A}_{p'}$ contains the sets of 
arithmetic non-rank progressions of all ${\cal A}_{p'},~5\leq 
p'\leq p.$   

The number $S(p)$ counts the non-ranks of ${\cal S}_p$ over one  
period $L(p).$ 
  
{\bf Definition~2.5.} Since~\cite{adhjw} $L(p)>S(p),$ there is 
a set ${\cal R}_p$ of {\it remnants} $r$ in its first period 
such that $r\not\in{\cal S}_p;$ they are twin-ranks or 
non-ranks to primes $p_j<p,$ where $p_j\ge 5$ is the $j$th 
prime. Let $M(j+1)=\frac{1}{6}(p^2_{j+1}-1).$ When all non-ranks 
to primes $p\le p_j$ are removed from the first period $[1,L(p_j)]$, 
all $r\le M(j+1)$ are twin ranks. These {\it front twin ranks} play 
the role of the primes $p\le\sqrt{x}$ in Eratosthenes' sieve that 
are left over when multiples of primes are removed. The prime $p_j$ 
is the twin sieve analog of $\sqrt{x}$ there; $p_j$ and $L(p_j)$ 
correspond to the variable $z$ and $P_z=\prod_{p\le z}p,$ 
respectively, in more sophisticated sieves.         
       
\section{Reworking the Twin Prime Formula}
 
If $p_j$ is the $j$th prime, then we need      
\begin{eqnarray}
L(p_j)=\prod_{5\leq p\leq p_j}p,~x=L(p_j)-M(j+1)  
\label{n1}
\end{eqnarray}
for the main result of Ref.~\cite{adhjw}, which 
is a Legendre-type formula for the number of 
twin ranks in the first period of length $L(p_j)$ 
of the supergroup ${\cal S}_{p_j},$ where $\pi_2$ 
counts twin pairs below $6x+1:$   
\begin{eqnarray}
\pi_2(6x+1)=R_0+\sum_{n\leq x,n|L_j(x)}
\mu(n)2^{\nu(n)}\bigg[\frac{x}{n}\bigg]+O(1)   
\label{leg}
\end{eqnarray}
where $[x/n]$ is the greatest integer function, 
$L_j(x)=\prod_{p_j<p\leq x}p,$ and  
\begin{eqnarray}
R_0=L(p_j)\prod_{5\leq p\leq p_j}(1-\frac{2}{p})
\sim \frac{Cx}{(\log\log x)^2},~C>0\cite{adhjw},~
p_j\sim \log x\to\infty 
\label{r0}
\end{eqnarray}
counts the number of remnants in ${\cal S}_{p_j},$ 
that is, twin ranks (prime pairs at distance $2$) 
{\bf and} non-ranks to primes $p_j<p\leq x.$ 
Therefore, the $n$ in the $\sum_n$ of Eq.~(\ref{leg}) 
run over these primes only and their products, and 
the upper limit is $x$ because the greatest integer 
function $[\frac{x}{n}]=0$ for $n>x.$ The twin pair 
counting function $\pi_2(M(j+1))$ is the number of 
front twin-ranks and the analog of $\pi(\sqrt{x})$ 
in Legendre's formula for the prime counting 
function $\pi(x)$ (see, e.g., pp. 2-3, Ch.~1 of 
Ref.~\cite{fi}); they are included in $R_0$. The 
error term $O(1)$ in Eq.~(\ref{leg}) accounts 
for the less than perfect cancellation at low 
values of $x$ of $R_0$ and the sum in 
Eq.~(\ref{leg}), but Eq.~(\ref{leg}) is only 
relevant at large $x$ in the following. Let us 
briefly sketch the cancellation of the too large 
$R_0$ against the sum in Eq.~(\ref{leg}) at large 
$\log\log x,$ upon decomposing $[x/n]=x/n-\{x/n\}$ 
as usual. Expanding the $R_0$ product into a sum 
and combining it with the corresponding sum of the 
ratios $x/n$ shifts the upper limit of the $R_0$ 
sum from $p_j$ to $x,$ when rewritten in its 
product form. This transforms its entire 
asymptotics from $\log\log x$ to $\log x$. For 
more details we refer to Theors.~5.7, 5.8 of 
Ref.~\cite{adhjw}

The asymptotic relations 
(\ref{r0}) derive from 
\begin{eqnarray}
\log L(p_j)=\sum_{5\leq p\leq p_j}\log p
=p_j+R(p_j)=\log x+O(\frac{\log^2 x}{x}), 
\end{eqnarray}   
where the error term comes from $M(j+1),$ and 
$R(p_j)$ is the remainder of the prime number 
theorem. 

The Dirichlet series characteristic of twin 
primes and associated with $R_0$ are 
\begin{eqnarray}
P_j(s)=\prod_{p>p_j}(1-\frac{2}{p^s})=
\prod_{p\leq p_j}(1-\frac{2}{p^s})^{-1}
\sum_{n=1}^{\infty}\mu(n)2^{\nu(n)}n^{-s}.  
\end{eqnarray}
They converge absolutely for $\sigma>1,$ as is 
evident from the majorant~\cite{ti} 
\begin{eqnarray}
\frac{\zeta^2(s)}{\zeta(2s)}=\sum_{n=1}^\infty
\frac{2^{\nu(n)}}{n^s},~\sigma>1.
\end{eqnarray} 
Note that $2^{\nu(n)}\sim \log n/\zeta(2)$ in 
the interval $[1,x]$ on average, which is 
shown in 4.4.18 of Ref.~\cite{rm}. The 
corresponding Dirichlet series for primes is 
\begin{eqnarray}
P_0(s)=\prod_{p\geq 2}(1-\frac{1}{p^s})=\frac{1}
{\zeta(s)}, 
\end{eqnarray}  
where $\zeta(s)=\sum_{n=1}^\infty\frac{1}{n^s}$ 
is the Riemann zeta function, and the analog of 
$R_0$ here is $x\prod_{p\le\sqrt{x}}(1-\frac{1}{p})$ 
there. 

We now use the Perron formula in essentially the 
form proved in 4.4.15 of Ref.~\cite{rm}. 

{\bf Lemma~3.1.} {\it Let the Dirichlet series 
$A(s)=\sum_{n=1}^\infty\frac{a_n}{n^s}$ be 
absolutely convergent for $\sigma=\Re(s)>1.$ 
Then}  
\begin{eqnarray}
\sum_{n\leq x} a_n=\frac{1}{2\pi i}\int_{
\sigma-iT}^{\sigma+iT}A(s)\frac{x^s}{s}ds+O\left(
\sum_{n=1,n\neq x}^\infty\left(\frac{x}{n}
\right)^\sigma |a_n|{\rm ~min}\left(1,\frac{1}
{T|\log\frac{x}{n}|}\right)\right),
\label{per} 
\end{eqnarray}
{\it where the lhs $\sum_{n\leq x}$ means that for 
$n=x,~a_n$ is reduced by $1/2.$}

{\bf Corollary~3.2.} {\it For} $\sigma>1$ 
\begin{eqnarray}\nonumber
&&\sum_{n\leq x} a_n\bigg[\frac{x}{n}\bigg]
=\frac{1}{2\pi i}\int_{\sigma-iT}^{\sigma+iT}A(s)
\zeta(s)\frac{x^s}{s}ds\\&&+O\left(\sum_{n=1,n\neq 
x}^\infty\left(\frac{x}{n}\right)^\sigma \left(
\sum_{d|n}|a_d|\right){\rm ~min}\left(1,\frac{1}
{T|\log\frac{x}{n}|}\right)\right).
\end{eqnarray}  

{\bf Proof.} This follows from Lemma~3.1 and the 
proof of 4.4.15 in Ref.~\cite{rm} using 
\begin{eqnarray}
\sum_{N\leq x}\sum_{n|N}a_n=\sum_{n\leq x}a_n
\bigg[\frac{x}{n}\bigg],~A(s)\zeta(s)=\sum_{
N=1}^\infty\frac{1}{N^s}\sum_{n|N}a_n.~\diamond
\end{eqnarray}

{\bf Lemma~3.3.}  
\begin{eqnarray}
P_1(s)\zeta^2(s)=(1-\frac{1}{2^s})^{-2}\prod_{p>2}
\left(1+\frac{1}{p^s(p^s-2)}\right)^{-1}=\frac{(1
-2^{-s})^{-2}}{D(s)}
\label{p1}
\end{eqnarray}
\begin{eqnarray}
D(s)=\prod_{p>2}\left(1+\sum_{\nu=0}^\infty\frac{
2^\nu}{p^{(\nu+2)s}}\right)=1+\sum_{N=4}^{\infty}
\frac{2^{2r_e(N)+2r_o(N)-2{\bar r}_e(N)-2{\bar 
r}_o(N)}}{N^s}
\label{hilf}
\end{eqnarray}
{\it converges absolutely for $\sigma>1/2.$ Here} 
\begin{eqnarray}
r_e(N)=\sum_{i=1}^m\nu_i,~r_o(N)=\sum_{i=1}^n
(\mu_i+3),~{\bar r}_e(N)=\sum_{\nu_i>0}1,~
{\bar r}_o(N)=\sum_{\mu_i>0}1
\end{eqnarray} 
{\it are additive functions for}
\begin{eqnarray}
N=p_{e_1}^{2(\nu_1+1)}\cdots p_{e_m}^{2(\nu_m+1)}
p_{o_1}^{2\mu_1+3}\cdots p_{o_n}^{2\mu_n+3}  
\label{enn}
\end{eqnarray}
{\it in Eq.~(\ref{hilf}).}

{\bf Proof.} Substituting in 
\begin{eqnarray}
\prod_{p>2}\frac{(1-\frac{2}{p^s})}{(1-\frac{1}
{p^s})^2}=\frac{1}{\prod_{p>2}(1+\frac{1}
{p^s(p^s-2)})}
\end{eqnarray} 
the expansions 
\begin{eqnarray}
\frac{1}{1-\frac{2}{p^s}}=1+\frac{2}{p^s}
+\frac{2^2}{p^{2s}}+\cdots, 
\end{eqnarray}
\begin{eqnarray}
1+\frac{1}{p^{2s}(1-\frac{2}{p^s})}=1+\sum_{
\nu=0}^\infty\frac{2^\nu}{p^{(\nu+2)s}}, 
\end{eqnarray}
yields Eq.~(\ref{hilf}) with $N$ of the form in 
Eq.~(\ref{enn}).~$\diamond$   

Thus for $\sigma>1$
\begin{eqnarray}\nonumber
P_j(s)^{-1}&=&\zeta^2(s)(1-\frac{1}{2^s})^2\prod_{
2<p\leq p_j}(1-\frac{2}{p^s})\prod_{p>2}\left(1+
\frac{1}{p^s(p^s-2)}\right)\\&=&\left(\frac{P_1(s)}
{\prod_{2<p\leq p_j}(1-\frac{2}{p^s})}\right)^{-1}
\label{pj}
\end{eqnarray}
with $P_1(s)$ from Eq.~(\ref{p1}). 

We now apply Cor.~3.2 to $P_j(s).$ This yields the 
Legendre-type formula {\bf before} it is split into 
its main and error terms according to Ref.~\cite{adhjw} 
so that the leading asymptotic term is $R_0$. 
 
{\bf Theorem~3.4.} {\it For} $\sigma>1,~R_0=L(p_j)
\prod_{5\leq p\leq p_j}(1-\frac{2}{p})$ and $x>0$ 
from Eq.~(\ref{n1}),  
\begin{eqnarray}\nonumber
&&\pi_2(6x+1)=R_0+\frac{1}{2\pi i}\int_{
\sigma-iT}^{\sigma+iT}\frac{(1-\frac{1}{2^s})^{-2}
x^sds}{s\zeta(s)\prod_{2<p\leq p_j}(1-\frac{2}
{p^s})D(s)}\\&&+O\left(\frac{\zeta^3(\sigma)
x^{\sigma}}{T}\right)+O\left(\frac{x\log^3 x}
{T}\right)+O(1),  
\label{we}
\end{eqnarray}
{\it with $D(s)$ from Eq.~(\ref{hilf}) and 
$T>0$ at least of order $x^c,~0<c<1.$}  
 
{\bf Proof.} We replace in Eq.~(\ref{leg}) the 
sum by the Perron integral of Cor.~3.2 with 
$A(s)=P_j(s)$ using Lemma~3.3 for $P_1(s)$ in 
conjunction with Eq.~(\ref{pj}). Canceling the 
factor $\zeta(s),$ this yields the Perron 
integral in Eq.~(\ref{we}).   

The Euler product of $D(s)$ in Eq.~(\ref{p1}) 
guarantees no zeros for $\sigma>1/2.$ Note that  
\begin{eqnarray}
\sum_{f|n}|\mu(f)|2^{\nu(f)}=\sum_{f|{\tilde n}}
2^{\nu(f)}=d_3({\tilde n})\le d_3(n), 
\end{eqnarray}  
where ${\tilde n}$ is the product of different 
prime divisors of $n$ and, for any $f|{\tilde n},$  
\begin{eqnarray}
d(f)=\sum_{\delta|f} 1=2^{\nu(f)}
\end{eqnarray} 
is the divisor function. Thus, we can use the 
majorant $d_3(n)$ in the error term in Cor.~3.2, 
where $\zeta^3(s)=\sum_{n=1}^{\infty}d_3(n)n^{-s}.$ 
We split the sum into three pieces as usual (see, 
e.g., Theor.~4.2.9 of Ref.~\cite{rm}) with $S_1=
\sum_{n<x/e}, S_2=\sum_{x/e<n<ex}, S_3=\sum_{n>ex}.$ 
For $S_1, S_3$ we have $|\log(x/n)|\ge 1.$ The 
total contribution due to $S_1$ and $S_3$ is at 
most $\zeta^3(\sigma)x^{\sigma}/T,$ which is the 
first error term in Eq.~(\ref{we}) with the 
constant $1$ implied by the $O(\cdots).$ 

For $S_2,$ we divide the sum into intervals of 
the type $I_k=[x\pm 2^kx/T,x\pm 2^{k+1}x/T]$ 
with $2^{k+1}/T<ex,$ and a shorter interval at 
the end if needed. The number of such intervals 
is $O(\log T).$ The contribution of the sum over 
such an interval to the remainder of Perron's 
formula is at most of order 
\begin{eqnarray}
\frac{1}{T}\sum_{I_k}d_3(n)\frac{T}{2^k}=
\frac{\sum_{I_k}d_3(n)}{2^k}.
\end{eqnarray}    
The length of $I_k$ is of order $2^kx/T,$ 
which is larger than $x^{1-c},$ since $T$ 
is at least of order $x^c$ for some $0<c<1.$  

Now recall the estimate (see Ref.~\cite{ti}, 
Ch.~12, formula~12.1.4):   
\begin{eqnarray}
\sum_{n<y}d_3(n)=yP_2(\log y)+O(y^{2/3}\log y),
\end{eqnarray}
$P_2$ being a certain polynomial of degree $2.$
 
It follows that 
\begin{eqnarray}
\sum_{I_k}d_3(n)=O(\frac{2^kx\log^2 x}{T})
+O(x^{2/3}\log x).
\end{eqnarray}
If we sum over $k$ the contribution of $S_2$ 
is at most of order $O(x\frac{\log^3 x}{T}),$ 
which gives the second error term in Eq.~(\ref{we}). 
The interval $(x/e<n<x)$ can be subdivided and 
treated similarly leading to the same bound. 
This completes the proof.~$\diamond$  
  
{\bf Corollary~3.5.} {\it The Riemann hypothesis 
(RH) is incompatible with the twin prime 
formula~(\ref{we}) of Theor.~3.4.}  

{\bf Proof.} Assuming RH, we shift the line of 
integration in Eq.~(\ref{we}) from $\sigma>1$ 
to $\sigma=\frac{1}{2}+\varepsilon$ for any 
$\varepsilon>0$ using Cauchy's theorem. Since 
RH implies the Lindel\"of hypothesis~\cite{ti}, 
we know that   
\begin{eqnarray}
\frac{1}{\zeta(s)}=O(|t|^\delta),~s=\sigma+it,~
\sigma\geq \frac{1}{2}+\varepsilon
\end{eqnarray}  
for some small $\delta>0$ that may depend on 
$\varepsilon.$ 
We note that the zeros $s_p=\log 2/\log p$ 
of $\prod_{p\le p_j}(1-2/p^s)$ cancel the 
corresponding poles of $D(s).$ Since only  
$s_3\approx 0.6309>0.5,$ we estimate for 
$\sigma\ge 1+\varepsilon$ 
\begin{eqnarray}
|(1-\frac{2}{3^s})D(s)|^{-1}=|(3^s-2+3^{-s})
\prod_{p\ge 5}[1+p^{-s}(p^s-2)^{-1}]|^{-1}
=O(1). 
\end{eqnarray} 
As $5\le p\le p_j\sim\log x$ in the product 
$\prod_p(1-2/p^s),$ the latter will be at 
most of order  
\begin{eqnarray}
\bigg|\prod_{5\le p\le p_j}(1-\frac{2}{p^s})\bigg|^{-1}
\le \prod_{5\le p\le p_j}(1-\frac{2}{\sqrt{p}})^{-1}
=O(\log x) 
\end{eqnarray}
for $\sigma\ge 1/2+\varepsilon$ as $p_j\sim\log 
x\to\infty.$ Hence, on $\sigma=\frac{1}{2}+
\varepsilon$ the vertical part of the Perron 
integral in Theor.~3.4 obeys 
\begin{eqnarray}
\int_{\frac{1}{2}+\varepsilon-iT}^{\frac{1}{2}+
\varepsilon+iT}\frac{(1-\frac{1}{2^s})^{-2}x^sds}
{s\zeta(s)\prod_{2<p\leq p_j}(1-\frac{2}{p^s})D(s)}
=O(T^\delta x^{\frac{1}{2}+\varepsilon}
\log T \log x), 
\label{rh}
\end{eqnarray} 
with the $\log T$ factor from the integration. 

On the horizontal line segments from 
$\frac{1}{2}+\varepsilon\pm iT$ to $\sigma\pm 
iT$ the Perron integral is bounded by 
$O(T^{\delta-1} x^\sigma).$ The factor $1/\log x$ 
from the integration cancels $\log x.$ 

The error terms of Theor.~3.4 are slightly 
smaller than these, respectively, and can be 
combined with them. Taking $\sigma=1+
\varepsilon,~T=x^\alpha$ and equating the 
exponents of $x$ in both error terms 
determines $\alpha=\frac{1}{2}.$ Therefore, 
the Perron integral plus error terms in 
Eq.~(\ref{we}) are of order $O(x^{\varepsilon+
(1+\delta)/2})$ and cannot reduce $R_0\sim C x/
(\log\log x)^2$ to the known bound~\cite{fi} 
$O(x/(\log x)^2)$ for $\pi_2(6x+1),$ q.e.a.
~$\diamond$ 

We next address the remainder of the twin prime 
formula~(\ref{leg}) after extracting its asymptotic 
law~\cite{adhjw} using the following Perron integral. 

{\bf Corollary~3.6.} {\it Let $A(s)$ be absolutely 
convergent for $\sigma>1,$ then for} $\sigma>1$   
\begin{eqnarray}\nonumber
&&\sum_{n<x}a_n\bigg\{\frac{x}{n}\bigg\}=\frac{1}
{2\pi i}\int_{\sigma-iT}^{\sigma+iT}ds x^s A(s)
\bigg[\frac{1}{s-1}-\frac{\zeta(s)}{s}\bigg]\\&&+
O\left(\sum_{n=1,n\ne x}\left(\frac{x}{n}
\right)^{\sigma}\sum_{d|n}|a_d|~{\rm min}
\left(1,\frac{1}{T|\log\frac{x}{n}|}\right)\right).  
\label{tpc}
\end{eqnarray}
 
{\bf Proof.} Using 
\begin{eqnarray}
\bigg\{\frac{x}{n}\bigg\}=\frac{x}{n}-\bigg[
\frac{x}{n}\bigg]
\end{eqnarray} 
and applying Lemma~3.1 to $xA(s+1)$ for the ratio 
$x/n,$ integrated along the line $\sigma>0,$ and 
Cor.~3.2 we obtain the Perron integral in 
Eq.~(\ref{tpc}) upon shifting $s\to s-1$ in the 
first term. Using $|a_n|\le \sum_{d|n}|a_d|,$ 
the error term of Lemma~3.1 combines with that of 
Cor.~3.2 giving that of Eq.~(\ref{tpc}).~$\diamond$

We now apply Cor.~3.6 to $P_j(s)$ which yields the 
Perron integral for the error term $R_E$ of 
Ref.~\cite{adhjw} after separating formula~(\ref{leg}) 
into its main and error terms so that the main term 
obeys the proper asymptotic law expected for twin 
primes~\cite{adhjw}. The error term is the same as 
in Theor.~3.4. For the cancellation involved in 
getting the proper asymptotics we refer to the 
discussion below Eq.~(\ref{r0}). Clearly, the 
sum in Eq.~(\ref{leg}), represented by the 
Perron integral in Theor.~3.4, is $-R_0$ plus 
an asymptotic term $cx/(\log x)^2$, with $c>0$ 
calculated in Ref.~\cite{adhjw}. Thus it is 
large, and an application of the contour 
deformation to the Perron integral in 
Theor.~3.4 into the known zero-free region of 
the Riemann zeta function fails to give a small 
value {\it unconditionally} because the optimal 
$a=0$ (as in the proof of Theor.~3.7 below) 
cannot be reached.   
       
{\bf Theorem~3.7.} {\it There are constants $a>0, 
0<b<c,1<\alpha<2$ so that the twin prime remainder 
takes on the form}
\begin{eqnarray}\nonumber
&&-R_E=\sum_{p_j<n<x,n|L_j(x)}\mu(n)2^{\nu(n)}
\bigg\{\frac{x}{n}\bigg\}+O(1)=O\left(\frac{
x^{1+\frac{a}{\log T}}\log^3 T}{T}\right)\\
\nonumber&&+O\left(\frac{x\log^3 x}{T}\right)
+O\left(x^{1-\frac{b}{\log T}}(\log T)^3
(\log x)^{\alpha}\right)+O\left(x^{1+\frac{a}
{\log T}}\frac{(\log T)^2(\log x)^{\alpha}}{T}
\right)\\&&=O\left(x\exp(-\sqrt{c\log x})
(\log x)^3\right),~T=\exp(\sqrt{c\log x}).  
\label{tpe}
\end{eqnarray}
{\bf Proof.} We start from Cor.~3.6 for $P_j(s)$ 
in conjunction with the error terms of Theor.~3.4:   
\begin{eqnarray}\nonumber
&&-R_E=\sum_{p_j<n<x,n|L_j(x)}\mu(n)2^{\nu(n)}
\bigg\{\frac{x}{n}\bigg\}+O(1)\\\nonumber&=&
\frac{1}{2\pi i}\int_{\sigma-iT}^{\sigma+iT}\frac{
(1-\frac{1}{2^s})^{-2}x^sds}{\zeta(s)\prod_{2<p\leq 
p_j}(1-\frac{2}{p^s})D(s)}\bigg[\frac{1}{(s-1)
\zeta(s)}-\frac{1}{s}\bigg]\\&+&O\left(\frac{
\zeta^3(\sigma)x^\sigma}{T}\right)+O\left(x\frac{
\log^3 x}{T}\right).
\end{eqnarray} 
By Chapt.~3, formula 3.11.8 of Ref.~\cite{ti} (this 
also follows from the sharper estimates in Lemma~12.3 
of Ref.~\cite{iv}) there is an absolute constant 
$c>0$ so that 
\begin{eqnarray}
\frac{1}{|\zeta(s)|}=O(\log(|t|+2),~0<t_0\leq |t|,~
\delta_t\leq\sigma\leq 1+\varepsilon,~\delta_t=1-
\frac{c}{\log(|t|+2)} 
\label{zbd}
\end{eqnarray} 
and $\frac{1}{|\zeta(s)|}=O(1)$ for $|t|<t_0,~
\sigma\geq\delta_t.$ Let $R_T$ be the rectangle 
joining the vertices  
\begin{eqnarray}
1+\frac{a}{\log T}-iT,~1+\frac{a}{\log T}+
iT,~\delta+iT,~\delta-iT,~\delta=1-\frac{c}{\log T}. 
\end{eqnarray}
We move the line segment of integration from $\sigma=1+
a/\log T,~a>0$ to the left on the line $\sigma=1-b/\log T$ 
with $a>0,~0<b<c$ to be chosen later. Then the bounds of 
$|\zeta(s)|^{-1}$ in Eq.~(\ref{zbd}) and below hold on 
the boundary of $R_T.$ The integrand is holomorphic 
inside and on the rectangle because $\zeta(s)$ does not 
vanish there and on the vertical line on the left it is 
of order at most 
\begin{eqnarray}
x^{1-\frac{b}{\log T}}\frac{\log T}{|s|\prod_{p\le p_j}
|1-p^{-1+\frac{b}{\log T}}|}.
\end{eqnarray}  
Since $p\le p_j\sim\log x$ we know that $p^{b/\log T}
\to 0$ as $x\to\infty$ provided $T$ grows with $x$ 
faster than any power of $\log x,$ which will be the 
case. Then a lower bound for the product will be at 
least of order $1/(\log p_j)^\alpha$ for any 
$1<\alpha<2.$ So the product is at most of order   
\begin{eqnarray}
|\prod_{2<p\le p_j}(1-\frac{2}{p^s})|^{-1}=
O((\log p_j)^{\alpha})=O((\log x)^\alpha). 
\end{eqnarray}
Integration over $s$ gives a factor $\log T.$  
Thus, the integral over the vertical segment 
is at most of order $O(x^{1-b/\log T}(\log 
T)^3(\log x)^\alpha).$ Similarly, the integrals 
over the horizontal segments are at most of order 
\begin{eqnarray}
O\left(x^{1+\frac{a}{\log T}}\frac{
(\log T)^2(\log x)^{\alpha}}{T}\right)    
\end{eqnarray}
Putting all this together we obtain the middle 
section of Eq.~(\ref{tpe}). 

Now we choose $T=\exp(\tau\sqrt{\log x})$ and 
optimize with respect to $a, b, \tau$ under the 
conditions $a>0, 0<b<c.$ In the limit we can set 
$b=c, a=0, \tau=\sqrt{c}$ and conclude with the 
bound on the rhs of Eq.~(\ref{tpe}). A 
comment on the choice $a=0$ is in order.  
The extra convergence factor $s/((s-1)\zeta(s))-1\to 
1/\zeta(s)-1$ at large $|t|$ is the Dirichlet series 
$\sum_{n\ge 2}\mu(n)n^{-s}$ that converges on 
$\sigma=1$ and oscillates rapidly at large $|t|$. 
This is how the real reduction from $[x/n]$ in 
Theor.~3.4 to $\{x/n\}$ in Theor.~3.7 plays out 
analytically. Its presence allows reaching the 
optimization point $a=0$ of the error terms 
representing the Perron integral in Theor.~3.7. 
This completes the proof.$~\diamond$

This proves that the (minimal) asymptotic law 
obtained in Ref.~\cite{adhjw} is valid with the 
remainder smaller than it by any positive power 
of $\log x.$  

\section{Summary and Discussion} 

When the Legendre-type formula for $\pi_2$ is 
reworked into a Perron integral involving 
$\zeta^{-1}(s),$ the nontrivial zeta zeros are 
seen to be linked to the twin prime counting 
function $\pi_2$. The asymptotic law of its 
leading term 
\begin{eqnarray}
R_0=L(p_j)\prod_{5\leq p\leq p_j}(1-\frac{2}
{p})\sim \frac{Cx}{(\log\log x)^2},~\log\log 
x\to\infty  
\end{eqnarray}
with $x=L(p_j)-M(j+1),~p_j\sim\log L(p_j)\sim\log x$ 
requires $\log\log x$ to become large. In contrast, 
only $\log x$ is large in the prime number 
theorem~\cite{rm},\cite{hw}. Therefore, the true 
asymptotic region of twin primes starts much higher 
up than for primes. Present numerical results of 
nontrivial zeta zeros have not yet reached the 
asymptotic twin prime realm. This is valid whether 
or not there are infinitely many twin primes, 
because $R_0$ is the number of remnants 
including twin pairs (i.e. twin ranks) {\bf and} 
non-ranks to primes $p_j<p<x.$ The Perron 
integral represents the latter's contributions 
that will reduce $R_0$ to $\pi_2(6x+1),~R_0$ 
being much larger than known bounds from sieve 
theory~\cite{rm},\cite{fi} on $\pi_2$ that are 
due to V. Brun, A. Selberg and others. 
Only nontrivial zeta zeros in the Perron integral 
can produce terms that reduce $R_0$ to the proper 
magnitude. Our first result is that the zeros on 
the critical line cannot do the job. Despite 
trillions of initial zeros on the critical 
line that are relevant for the prime number 
distribution without asymptotic twin prime 
attributes, once twin prime asymptotics matter 
zeta zeros must move off the critical line  
toward the borders of the critical strip.                    

Finally, from the point of view of our twin prime 
formulas~(\ref{leg}),(\ref{we}) a finite number of 
twin primes is neither a simple nor natural case, 
as it would require the cancellation of the leading 
and all subleading asymptotic terms involving 
fine-tuning of the large primes ($>p_j$) that 
organize the non-ranks. 

But this never happens because, when the Perron 
integral is developed for $\sum_{n<x}\mu(n)2^{
\nu(n)}\{x/n\}$ using $P_j(s)$ in the known 
zero-free region of the Riemann zeta function, 
the twin prime theorem near primorial arguments 
follows, our second result.    

Third, this analysis extends--mutatis mutandem--
to all other twin prime cases 
\cite{hjw1},\cite{hjw2},\cite{hjw3} of the 
classes I, II, III of the 
classification~\cite{hjw},\cite{hjwe}, that is 
to say that any twin prime case prevents RH from 
being valid, except for the initial, extremely 
long stretch, and has a (minimal) asymptotic 
law of the expected form at primorial arguments.

%\end{references}

\begin{thebibliography}{}

\bibitem{adhjw} H. J. Weber, Twin Prime Sieve, 
www.arxiv.org/1203.5240. 

\bibitem{hri} H. Riesel, Prime Numbers and 
Computer Methods for Factorization, 2nd ed., 
Birkh\"auser, Boston, 1994. 

\bibitem{rm} M. Ram Murty, Problems in 
Analytic Number Theory, Springer, New York, 
2001. 

\bibitem{hw} G. H. Hardy and E. M. Wright,  
An Introduction to the Theory of Numbers, 
5th ed., Clarendon Press, Oxford, 1988.

\bibitem{fi} R. Friedlander and H. Iwaniec, 
Opera Cribro, Amer. Math. Soc. Colloq. Publ. 
57, Prov. RI, 2010. 

\bibitem{ti} E. C. Titchmarsh, The Theory of the 
Riemann Zeta-Function, 2nd ed. edited by D. R. 
Heath-Brown, Oxford, Clarendon Press, 1986. 

\bibitem{iv} A. Ivi\'c, The Riemann 
Zeta-Function, Dover, Mineola, N.Y., 2003. 

\bibitem{hjw1} H. J. Weber, Sieves for Twin
Primes in Class I, www.arxiv.org/1204. 5728.

\bibitem{hjw2} H. J. Weber, Sieves for Twin
Primes in Class II, UVa preprint, 2012.

\bibitem{hjw3} H. J. Weber, Sieves for Twin
Primes in Class III, UVa preprint, 2012.

\bibitem{hjw} H. J. Weber, Regularities of
Prime Number Twins, Triplets and Multiplets,
Global. J. Pure Applied Math., in press, www. 
adsabs. harvard. edu/abs/2011 arXiv1103.0447W. 

\bibitem{hjwe} H. J. Weber, Exceptional 
Prime Number Twins, Triplets and Multiplets,
www.arxiv.org/1102.3075. 

\end{thebibliography}
\end{document}